\documentclass[11pt]{amsart}
\pdfoutput=1
\usepackage[utf8]{inputenc}
\usepackage{amsmath, amssymb}
\usepackage[english]{babel}
\usepackage{graphicx}
\usepackage{amsthm}
\usepackage{array}
\usepackage{amsthm}
\usepackage{mathtools}
\usepackage{thmtools}
\usepackage{enumitem}
\usepackage[margin=.85 in]{geometry}
\allowdisplaybreaks
\usepackage{lpic}

\usepackage[dvipsnames]{xcolor}
\usepackage{hyperref}
\definecolor{candyapplered}{rgb}{1.0, 0.03, 0.0}
\definecolor{mediumblue}{rgb}{0.0, 0.0, 0.8}

\hypersetup{
bookmarksnumbered,
unicode=false,          
pdftoolbar=true,        
pdfmenubar=true,        
pdffitwindow=false,     
pdfstartview={FitH},    
pdfauthor={Jonah Gaster},
pdftitle={Boundary slopes in the Markov ordering},
bookmarksnumbered,
colorlinks=true,
linkcolor=mediumblue,
citecolor=candyapplered,
urlcolor=blue}

\usepackage{caption,subcaption}
\usepackage[ocgcolorlinks]{ocgx2}

\declaretheorem[numberwithin=section]{theorem}

\declaretheorem[sibling=theorem]{lemma}
\declaretheorem[sibling=theorem, style=remark]{remark}
\declaretheorem[sibling=theorem]{proposition}

\declaretheorem[sibling=theorem, style=remark]{example}

\declaretheoremstyle[
spaceabove=6pt,
spacebelow=6pt,
headfont=\normalfont\bfseries,
notefont=\mdseries,
notebraces={(}{)},
bodyfont=\normalfont,
postheadspace=.5em,
numbered=no]{mystyle}
\declaretheorem[style=mystyle]{problem}

\numberwithin{equation}{section}

\newcommand{\bZ}{\mathbb{Z}}

\newcommand{\bN}{\mathbb{N}}
\newcommand{\bH}{\mathbb{H}}

\newcommand{\bR}{\mathbb{R}}
\newcommand{\bQ}{\mathbb{Q}}

\newcommand{\SL}{\mathrm{SL}}
\newcommand{\PSL}{\mathrm{PSL}}
\newcommand{\cF}{\mathcal{F}}
\newcommand{\cH}{\mathcal{H}}
\newcommand{\cM}{\mathcal{M}}

\newcommand{\cQ}{\mathcal{Q}}

\newcommand{\cL}{\mathcal{L}}

\newcommand{\tr}{\mathrm{tr\,}}

\newcommand{\upd}{\mathrm{d}}

\mathtoolsset{showonlyrefs=true}

\begin{document}
\title[Boundary slopes for the Markov ordering]{Boundary slopes for the\\ Markov ordering on relatively prime pairs}
\author[Gaster]{Jonah Gaster}
\date{August 9, 2021}

\address{Department of Mathematical Sciences, University of Wisconsin-Milwaukee}
\email{gaster@uwm.edu}

\keywords{Curves on surfaces, hyperbolic geometry, Markov numbers}

\begin{abstract}
Following McShane, we employ the stable norm on the homology of the modular torus to investigate the Markov ordering on the set of relatively prime integer pairs $(q,p)$ with $q\ge p\ge0$.
Our main theorem is a characterization of slopes along which the Markov ordering is monotone with respect to $q$, confirming conjectures of Lee-Li-Rabideau-Schiffler that refine conjectures of Aigner.
The main tool is an explicit computation of the slopes at the corners of the stable norm ball for the modular torus. 
\end{abstract}

\maketitle

\section{Introduction}

The set of \emph{Markov numbers} $\cM=\{1,2,5,\ldots\}$ is the set of positive integers that appear in a Markov triple, i.e.~a solution $(x,y,z)\in\bN^3$ to the cubic
\begin{equation}
\label{eq:markov}
x^2+y^2+z^2=3xyz~.
\end{equation}
The \emph{Markov Uniquness Conjecture} (MUC), apparently due to Frobenius, asserts that each positive Markov triple is determined by its largest entry \cite{Aigner,CF,Frobenius,Markoff}.

The Markov numbers admit an interpretation as a system of labels on the set $\cH$ of horoball regions in the complement of a rooted planar trivalent tree.
That is, there is a surjective map $\lambda_{\cM}:\cH \to \cM$, whose injectivity is equivalent to MUC \cite[p.~54]{Aigner}. 
There is another well-known natural parametrization of $\cH$ via Farey fractions $\lambda_{\cF}:\cH \cong \bQ\cap[0,1]$, so one may obtain the Markov numbers as $m_{p/q}=\lambda_{\cM}\circ\lambda_{\cF}^{-1}(p/q)$. See \autoref{fig:lambdaM} and \autoref{fig:lambdaF} for an illustration of $\lambda_\cM$ and $\lambda_\cF$.

It will be convenient to record this map as follows: let $\cQ\subset \bZ^2$ indicate the set of relatively prime pairs $(q,p)$ with $q\ge p\ge 0$. 
The map $(q,p)\mapsto m_{p/q}$ records the Markov label associated to the horoball with label $p/q$. MUC predicts that this map is in fact a bijection; because $\cM\subset \bN$, one would obtain a somewhat mysterious ordering of $\cQ$.
We will refer to this ordering as the `Markov ordering', denoted $\prec_{\cM}$. (For the careful reader who prefers not to assume MUC, $\prec_{\cM}$ is a priori only a strict partial order, with potentially nontrivial mutually incomparable subsets for horoballs with common Markov labels.) 
See \autoref{fig:orderQ} for an illustration.

Initiating a study of $\prec_\cM$, Aigner offered some conjectural statements suggested by numerical experiment \cite{Aigner}: assuming the following pairs are relatively prime,

\begin{enumerate}[label=(\roman*), leftmargin=2.5cm]
\item $(q,p) \prec_\cM (q',p)$ when $q'>q$,
\item $(q,p) \prec_\cM (q,p')$ when $p'>p$, and
\item $(q,p) \prec_\cM (q',p')$ when $q'>q$ and $p+q=p'+q'$.
\end{enumerate}
These three statements were subsequently proven by several authors using the theory of cluster algebras and continued fractions \cite{Rabideau,RS,LLRS,LPTV}. 
Shortly thereafter, Greg McShane gave a unified proof of Aigner's conjectures \cite{McShane} by exploiting hyperbolic geometry and the so-called \emph{stable norm}, a norm on the homology of the modular torus induced by the hyperbolic length function \cite{MR1}.

\begin{figure}
\centering
\begin{minipage}[t][7.cm][t]{.33\textwidth}
\centering
\vspace{1cm}
\begin{lpic}{tree1(4.5cm)}
\large
\lbl[]{28,65;$1$}
\lbl[]{84,65;$2$}
\lbl[]{56,18;$5$}
\lbl[]{4,18;$13$}
\lbl[]{105,18;$29$}
\lbl[]{28,-8;$194$}
\lbl[]{84,-8;$433$}
\end{lpic}
\vspace{1cm}
\caption{The\\ labelling $\lambda_\cM$.}
\label{fig:lambdaM}
\end{minipage}\hfill
\begin{minipage}[t][7.cm][t]{.33\textwidth}
\centering
\vspace{1cm}
\hspace{.2cm}
\begin{lpic}{tree1(4.5cm)}
\large
\lbl[]{28,65;$0/1$}
\lbl[]{84,65;$1/1$}
\lbl[]{56,18;$1/2$}
\lbl[]{2,18;$1/3$}
\lbl[]{106,18;$2/3$}
\lbl[]{28,-8;$2/5$}
\lbl[]{84,-8;$3/5$}
\end{lpic}
\vspace{1cm}
\caption{The\\ labelling $\lambda_\cF$.}
\label{fig:lambdaF}
\end{minipage}\hfill
\begin{minipage}[t][7.cm][t]{.3\textwidth}
\centering
\vspace{1cm}
\begin{lpic}{orderQ(4.5cm)}
\end{lpic}
\caption{$\prec_\cM$ on $\cQ$.}
\label{fig:orderQ}
\end{minipage}\hfill
\end{figure}

Aigner's conjectures concern the monotonicity of $\prec_\cM$ along lines of slopes $\infty$, $0$, and $-1$.
This analysis was pushed further in the work of Lee-Li-Rabideau-Schiffler, who consider monotonicity of $\prec_\cM$ along lines of other slopes. They show \cite[Thm.\,1.2]{LLRS}:
\begin{enumerate}[label=(\roman*), leftmargin=2.5cm]
\item $\prec_\cM$ is monotone increasing with $q$ along lines of slope $\ge -\frac87$, 
\item$ \prec_\cM$ is monotone decreasing with $q$ along lines of slope $\le -\frac54$, and 
\item $\exists$ lines of slopes $-\frac76$ and $-\frac65$ along which $\prec_\cM$ is not monotone.
\end{enumerate}
Moreover, the authors make precise conjectures concerning optimality of slopes for monotonicity of $\prec_\cM$ \cite[Conj.\,6.8, 6.11, 6.12]{LLRS}\footnote{While an initial draft of \cite{LLRS} does not contain these conjectures, they can be found in a revised manuscript which was shared in private communication.
This author thanks those authors for their openness.}. 
The purpose of this note is to confirm these conjectures (though we stop just short of the full statement of \cite[Conj.\,6.12]{LLRS}, 
which reaches beyond $\cM$ to an order that concerns non-relatively prime pairs as well -- see \autoref{rem:Conj612}).

There is a well-known geometric interpretation of $\cM$, first discovered by Gorshkov and Cohn \cite{Gorshkov, Cohn, Series}, which we now recall:
Let $X$ be the modular torus.
The fundamental group of $X$ can be lifted to $\SL(2,\bZ)$ so that the set of traces of simple closed geodesics on $X$ is exactly $3\cM$.
Because traces and hyperbolic lengths are related by a monotone function \eqref{eq:length-Markov}, one finds that $\prec_\cM$ coincides with the order on $\cQ$ induced by the hyperbolic length function (that is, with the appropriate identification of $H_1(X,\bZ) \approx \bZ^2$, in which the three shortest geodesics on $X$ have representatives $(1,0), (0,1), (-1,1)$).

Recently, McShane has shown how to use this viewpoint to prove Aigner's conjectures, relying on a certain convexity induced by the hyperbolic length function.
The stable norm $\|{\cdot}\|_s$ is a norm on $\bR^2\approx H_1(X,\bR)$ whose value on primitive integral points coincides with the hyperbolic length of the corresponding simple geodesic (see \cite{McShane,MR1,SoV,BBI} for more detail).
McShane-Rivin observed that the boundary of the $\|{\cdot}\|_s$-ball is strictly convex \cite{MR1}, and Aigner's conjectures about $\prec_\cM$ can then be deduced rather quickly by consideration of the unit $\|{\cdot}\|_s$-ball \cite{McShane}.

We follow McShane's technique and Lee-Li-Rabideau-Schiffler's conjectures:

\begin{theorem}
\label{thm:main thm}
Let $r\in\bQ$, $\sigma_-=\displaystyle -\frac
{\log \left(\frac32 + \frac12\sqrt 5\right)}
{\log \left(\frac32+\frac3{10}\sqrt 5\right)}
\approx -1.2417$, and $\sigma_+=
\displaystyle -\frac
{\log \left(\frac32+\frac34\sqrt 2\right)}
{\log \left(\frac43+\frac23\sqrt 2 \right)}\approx -1.1432$.

\begin{enumerate}[label=\normalfont(\roman*), leftmargin=1.2cm]
\item
\label{it:sigma-}
If $r< \sigma_-$, then $\prec_\cM$ is monotone decreasing with $q$ along lines of slope $r$.
\item
\label{it:sigma+}
If 
$r> 
\sigma_+$, then $\prec_\cM$ is monotone increasing with $q$ along lines of slope $r$.
\item
\label{it:nope}
If $r\in (\sigma_-,\sigma_+)$, there exist lines of slope $r$ along which $\prec_\cM$ is not monotone.
\end{enumerate}
\end{theorem}

\begin{remark}
\cite[Thm.\,1.2]{LLRS} is an immediate corollary, since $-\frac54 < \sigma_- < -\frac65 < -\frac76 < \sigma_+ < -\frac87$.
\end{remark}

\begin{remark}
\label{rem:antimodal}
In fact, the proof of \autoref{thm:main thm}\ref{it:nope} demonstrates slightly more: Let $\cL$ be a line of slope $r\in (\sigma_-,\sigma_+)$. 
There exist arbitarily large $t\in \bR$ so that the scaling $\cL_t:=t\cdot \cL$ has intersection $\cL_t \cap \cQ$ which is nonempty and strictly antimodal with respect to $q$. 
That is, for arbitrarily large $t$ one has $\cL_t\cap \cQ= \{(q_1,p_1),\ldots,(q_n,p_n)\}$ with $n\ge 3$ and $q_1<\ldots<q_n$, and there exists $j\in\{2,\ldots, n-1\}$ so that
$(q_i,p_i)\prec_\cM (q_j,p_j)$ for $i<j$ and $(q_j,p_j)\prec_\cM (q_i,p_i)$ for $i>j$.
See \autoref{rem:strict antimodality} and \autoref{rem:Conj612}.
\end{remark}

The proof of \autoref{thm:main thm} amounts to a computation about (one-sided) slopes formed at two corners of the stable norm ball in the homology of the modular torus. 
This same computation can be performed at any point of rational slope (see \autoref{prop:corners} below). 
In fact, this allows one to make the following more precise statement about $\prec_\cM$:

\begin{theorem}
\label{thm:fixed point}
Let $(q,p),(q',p')\in\cQ$, and let $\mu^+_{p/q}$ and $\mu^-_{p/q}$ be given by \eqref{eq:plus} and \eqref{eq:minus}.

If $q'<q$ and $\displaystyle \frac{q'-q}{p'-p} \le \mu^+_{p/q}$ then $(q,p)\prec_\cM (q',p')$.

If $q'>q$ and $\displaystyle \frac{q'-q}{p'-p} \ge \mu^-_{p/q}$ then $(q,p) \prec_\cM (q',p')$.
\end{theorem}
\begin{remark}
The constants appearing in \autoref{thm:main thm} are given by $\sigma_-=\mu^+_{0/1}$ and $\sigma_+=\mu^-_{1/1}$.
\end{remark}

\begin{figure}
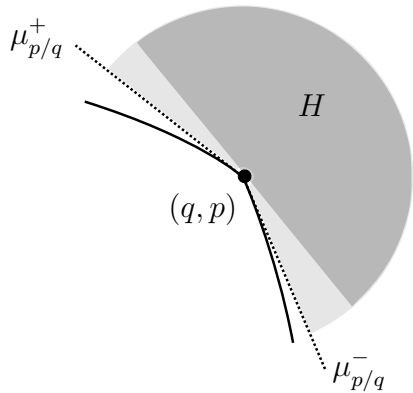

\centering
\begin{lpic}{corner(4.5cm)}
\large
\lbl[]{13.5,17;$(q,p)$}
\lbl[]{-4,35;$\mu^+_{p/q}$}
\lbl[]{30,0;$\mu^-_{p/q}$}
\lbl[]{25,28;$H$}
\end{lpic}
\caption{A support plane $H$ at a corner of the stable norm ball.}
\label{fig:corner}
\end{figure}

\begin{proof}
If either $q'<q$ and $\frac{q'-q}{p-'p}\le \mu_{p/q}^+$ or $q'>q$ and $\frac{q'-q}{p'-p}\ge \mu^-_{p/q}$, there is a half plane $H$, containing $(q',p')$, which is a support plane for the $\|{\cdot}\|_s$-ball of radius $\|(q,p)\|_s$ centered at the origin. Strict convexity implies that $(q',p')$ is therefore outside the ball, so $\|(q,p)\|_s\le \|(q',p')\|_s$.
See \autoref{fig:corner}.
\end{proof}

Our computation of the extremal slopes above is convenient when phrased in terms of Fock's function $\Psi$.
Vladimir Fock introduced the function $\Psi:[0,1/2]\to\bR$ whose value at rationals is given by 
\begin{equation}
\label{eq:psi}
\Psi(p/q) = \frac{\cosh^{-1}\left(\frac32 m_{T\left(\frac pq\right)}\right)}{q}~,
\end{equation}
where $T\left( \frac pq \right) = \frac p{q-p}$.
Suitably interpreted, $\Psi(p/q)$ is the ratio of (one-half of) the hyperbolic length to an `arithmetic height' for a simple closed geodesic corresponding to the homology class $T(q,p)=(q-p,p)$.
Fock showed that $\Psi$ is a continuous, convex, and decreasing function \cite{Fock}.

\begin{remark}
\label{rem:T}
We hope the reader will allow the slightly abusive repetition of $T$; evidently, under the natural bijection $\cQ\approx \bQ\cap[0,1]$ these two maps are `the same'.
\end{remark}

Recall that, if $A\in \SL(2,\bR)$ is a hyperbolic transformation with translation distance $\delta(A)$, one has
\[
2\cosh\frac{\delta(A)}2 = \tr(A)~.
\]
When $\gcd(p,q)=1$, the translation distance of the holonomy of $(q,p)$ is the length $\|(q,p)\|_s$ of its simple geodesic representative, and the trace is three times the corresponding Markov number $m_{p/q}$, so we find 
\begin{equation}
\label{eq:length-Markov}
\| (q,p) \|_s= 2\cosh^{-1}\left( \frac32 m_{p/q}\right)~.
\end{equation}
An equivalent description of the value of $\Psi$ at $p/q$ is thus given by 
\begin{equation}
\label{eq:psi-snorm}
\Psi\left( \frac pq \right) = \frac{ \| T(q,p) \|_s}{2q}~.
\end{equation}

\begin{remark}
\label{rem:Fock normalization}
The definition of $\Psi$ given in \eqref{eq:psi} and \eqref{eq:psi-snorm} looks cosmetically distinct from others in the literature (\cite[eq.\,(3)]{SoV} for example), owing to `$T(p/q)$' in place of `$p/q$' (or `$T(q,p)$' in place of `$(q,p)$'). 
This is the result of two distinct normalizations for the Farey labelling $\lambda_\cF$ in the literature. 
In one, the root of the tree is incident to $0/1$, $1/0$ and $1/1$, as in \autoref{fig:lambdaF} and \cite{Aigner,LLRS,McShane} (observe that this is natural by virtue of the integral homology $\bZ^2$), and in the other the root sees $0/1$, $1/1$, and $1/2$, as in \cite{Fock, SpV, SoV} (somewhat more natural when considering the holonomy of the modular torus in $\PSL(2,\bZ)$).
We've chosen a presentation of $\prec_\cM$ that follows \cite{Aigner,LLRS,McShane}, 
but our function $\Psi$ matches that of \cite{Fock, SoV, SpV}.
\end{remark}

The structure of $\Psi$ is closely tied to the Markov numbers and their cousins, the Lagrange numbers \cite{SpV,SoV}.
Sorrentino-Veselov asked for an explicit computation of the left- and right-derivatives of $\Psi$ at rationals \cite[\S5]{SoV}. 
We provide an answer below.

\begin{proposition}
\label{prop:psi derivatives}
Suppose that $\frac pq$, $\frac {r_1}{s_1}$, and $\frac{r_2}{s_2}$ span a Farey triangle with $\frac{r_1}{s_1}<\frac pq < \frac{r_2}{s_2}$, and let $n_i := m_{T(r_i/s_i)}$ and $n:=m_{T(p/q)}$, so that $(n,n_1,n_2)$ is a Markov triple.

The left- and right-derivatives of $\Psi$ at $p/q$ can be computed as follows:
\begin{equation}
\label{eq:leftPsi}
\left.\frac{\upd }{\upd t}\right|_{t=\frac pq^-} \Psi(t) = 
-s_2 
\cosh^{-1}\left( \frac{3n}2\right)
- q \log\left( \frac {3n_2}2-
\frac 
{9n n_2-6n_1}
{2\sqrt{9n^2-4}}
\right)
\end{equation}
\begin{equation}
\label{eq:rightPsi}
\left.\frac{\upd }{\upd t}\right|_{t=\frac pq^+} \Psi(t) = 
s_1 
\cosh^{-1}\left( \frac{3n}2\right)
+ q \log\left( \frac {3n_1}2-
\frac 
{9n n_1-6n_2}
{2\sqrt{9n^2-4}}
\right)
\end{equation}
\end{proposition}

\begin{remark}
It has come to the attention of the author that the calculation of the left- and right-derivatives of $\Psi$ was carried out recently, in much the same fashion, by Robert Hines \cite[p.\,6]{Hines}. 
\end{remark}

One can translate this computation about $\Psi$ into a computation of the slopes at the corners of spheres with respect to the stable norm: 

\begin{proposition}
\label{prop:corners}
Suppose $(q,p)\in\cQ$. Let $L=\left.\frac{\upd }{\upd t}\right|_{t=T^{-1}\left(\frac pq\right)^-} \Psi(t)$, $R=\left.\frac{\upd }{\upd t}\right|_{t= T^{-1}\left(\frac pq\right)^+} \Psi(t)$, whose values are recorded in \eqref{eq:leftPsi} and \eqref{eq:rightPsi}, and let $\ell=\frac12\|(q,p)\|_s$ (equivalently, $\ell=\cosh^{-1}(\frac32 m_{p/q})$).

The corner of the stable norm ball at $(q,p)$ has left slope $\mu^+_{p/q}$ and right slope $\mu^-_{p/q}$ given by
\begin{align}
\label{eq:plus}
\mu^+_{p/q} & = -\frac{\ell - Rp}{\ell + R q}, \text{ and }\\
\label{eq:minus}
\mu^-_{p/q} & = -\frac{\ell- Lp}{\ell +L q} ~.
\end{align}
\end{proposition}

\begin{example}
\label{ex:minus}
We compute $\mu^+_{0/1}$.
In this case, $\ell=\cosh^{-1}\frac 32=\log\left( \frac{3+\sqrt 5}2\right)$. In order to determine the right-derivative $R$ of $\Psi$ at $0/1$, we consider the Farey trio $\left(\frac{r_1}{s_1} , \frac pq,\frac{r_2}{s_2}\right)=\left(\frac10,\frac01,\frac11\right)$ and corresponding Markov triple $(n_1,n,n_2)=(1,1,1)$. By \eqref{eq:rightPsi} we have
\begin{align*}
R & = 0\cdot \cosh^{-1} \frac 32 +1\cdot \log\left( \frac32 - \frac{3}
{2\sqrt{5}}\right) 
= \log \left( \frac32-\frac3{10}\sqrt 5\right)~.
\end{align*}
Therefore by \autoref{prop:corners} we have
\begin{align*}
\mu^+_{0/1} = -\frac{\log\frac{3+\sqrt 5}2}{\log\frac{3+\sqrt 5}2+\log \frac{15-3\sqrt 5}{10}} = -\frac{ \log\left( \frac32+\frac12\sqrt 5\right)}{\log\left( \frac 32 +\frac3{10}\sqrt 5\right)}~.
\end{align*}
\end{example}

\begin{example}
\label{ex:plus}
We compute $\mu^-_{1/1}$. To compute the left-derivative of $\Psi$ at $1/2$ (note that $1/2=T^{-1}(1/1)$), fix Farey trio $\left(\frac{r_1}{s_1}, \frac pq,\frac{r_2}{s_2} \right)=\left(\frac01,\frac12,\frac11\right)$ and corresponding Markov triple $(1,2,1)=(n_1,n,n_2)$. By \eqref{eq:leftPsi},
\begin{align*}
L&=-1\cdot \cosh^{-1}3 -2\cdot \log \left( \frac32 - \frac{12}{8\sqrt{2}}\right) 
=\log \frac89~.
\end{align*}
Therefore by \autoref{prop:corners} we have
\begin{align*}
\mu^-_{1/1} = -\frac{\cosh^{-1}3 - \log \frac 89}{\cosh^{-1}3+\log \frac 89}= -\frac{\log\frac98(3+2\sqrt 2)}{\log\frac89(3+2\sqrt 2)}=-\frac{\log\frac3{2\sqrt 2}(1+\sqrt 2)}{\log\frac{2\sqrt 2}3(1+\sqrt 2)}=-\frac{\log\left(\frac32+\frac34\sqrt2\right)}{\log\left(\frac43+\frac23\sqrt 2\right)}~.
\end{align*}
\end{example}

\begin{remark}
McShane-Rivin offer a more general version of \autoref{prop:corners} \cite[Thm.\,2.1]{MR2}, though the proof has not appeared in the literature. Their theorem states that the exterior angle at $(q,p)$ decreases exponentially with $\max\{p,q\}$, a fact that can be deduced from the formulas above. Though \autoref{prop:corners} has the advantage of being exact, we emphasize that the McShane-Rivin theorem is more general: Their control on the boundary slopes applies to any hyperbolic punctured torus, and also accounts for the slope at irrational points on the boundary of the stable norm ball.
It might be interesting to find analogues of \eqref{eq:leftPsi} and \eqref{eq:rightPsi} at points on the boundary of the stable norm ball corresponding to irrational laminations (cf.~\cite[\S5]{SoV}). 
\end{remark}

We conclude this discussion with a challenge: observe that the number of mutually $\prec_\cM$-incomparable elements of $\cQ$ is a bound for the potential non-injectivity of $\lambda_\cM$ in relation to MUC (indeed MUC predicts that $\lambda_\cM$ is injective and that $\prec_\cM$ is a total order). 
Mutually $\|{\cdot}\|_s$-incomparable points in $\cQ$ are close to a subject that has attracted attention in convex analysis, as they can be thought of as rational points with common denominator on the $\|{\cdot}\|_s$-unit sphere, a planar convex curve. 
By \cite[Thm.~1]{Petrov}, the number of such elements can be bounded as
\[
\# \lambda_\cM^{-1}(n) = o\left( (\log n)^{2/3} \right)~.
\]
The latter we refer to as \emph{Petrov's bound}.

Notice that Petrov's bound is considerably stronger in the worst-case scenario than the super-logarithmic bound $\approx 2^{\frac{\log n}{\log \log n}} $ that arises from counting roots of $-1$ modulo $n$, as in \cite[Thm.\,2.19]{Aigner}, or the logarithmic bound $O(\log n)$ that follows from the analysis of Zagier \cite{Zagier} (see also \cite[Cor.\,1.5(b)]{LLRS} for another logarithmic bound).
An error-term for the count of simple closed geodesics on $X$ was obtained by McShane-Rivin in \cite{MR2} which provides $\#\lambda_\cM^{-1}(n)=O(\log n \,\log\log n)$, again super-logarithmic, although they conjecture a much stronger error term \cite[Conj.\,4.3]{MR2} which would imply the much stronger bound $O\left(\left(\log n\right)^{\frac12+\epsilon}\right)$, for any $\epsilon>0$. As far as the author is aware, this McShane-Rivin conjecture is far beyond current methods (e.g.\,\cite[Thm.\,1.1]{EMM}).

\begin{problem}
Can Petrov's bound be improved for the $\|{\cdot}\|_s$-unit sphere? 
Better bounds for strictly convex \emph{smooth} plane curves \cite{Swinnerton-Dyer} are evidently unapplicable here, but can one leverage McShane-Rivin's observation that it is \emph{infinitely flat} at points of irrational slope \cite[Thm.\,2.1]{MR2} to deduce improvements to Petrov's bound?
Are better bounds for $\#\lambda_\cM^{-1}(n)$ available by some other means?
\end{problem}

\begin{remark}
According to \cite{Plagne}, there exist strictly convex curves with rational points of denominator $q$ growing just slower than $q^{2/3}$ (for certain choices of $q$).
Therefore, any improvements to Petrov's bound must rely on more than merely strict convexity of the $\|{\cdot}\|_s$-ball.
\end{remark}

\subsection*{Acknowledgements}
The author thanks Greg McShane and Li Li for helpful correspondence.

\bigskip

\section{Hyperbolic geometry of Markov triples on the modular torus}
The modular torus $X$ (also called the \emph{equianharmonic torus} in the literature) can be obtained as the quotient of $\bH^2$ by the commutator subgroup of $\PSL(2,\bZ)$. The holonomy group of $X$ admits a lift to $\SL(2,\bZ)$ (in fact, such a lift exists in much greater generality \cite{Culler}), so one may speak of the `trace of a simple closed geodesic on $X$'.
For punctured tori, any such choice of lift will give rise to a representation sending the cusp to an element of trace $-2$, from which one can deduce using the $\SL_2$-trace identity that triples of simple closed curves pairwise intersecting once have traces that satisfy a Markov-like equation. 
A lift exists so that all simple closed curves have positive integral traces, and one finds that the set of traces of simple closed geodesics on $X$, in the chosen lift, will be precisely $3\cM$. See \cite{Goldman} for details.
We indicate the trace of the simple closed geodesic $\gamma$ on $X$ by $3m_\gamma$, for $m_\gamma\in \cM$.

In order to compute derivatives for $\Psi$, we recall some aspects of simple closed geodesics on $X$.
The following is a topological exercise.  
See \cite{McShane-Weierstrass} for detail.

\begin{lemma}
\label{lem:hyp geom}
Any collection of three simple closed geodesics pairwise intersecting once on $X$ have 
two complementary triangles, exchanged by the elliptic involution, whose vertices are the three Weierstrass points of $X$ and whose side lengths are the half-lengths of the three geodesics.
\end{lemma}

In order to compute derivatives of Fock's function, it will be important to keep track of the behavior of hyperbolic lengths under Dehn twisting. We indicate the mapping class obtained via the right Dehn twist about the simple closed curve $\gamma$ by $\tau_\gamma$, and we indicate the half-length of $\gamma$ on $X$ by $\ell(\gamma)$.

\begin{figure}[h]
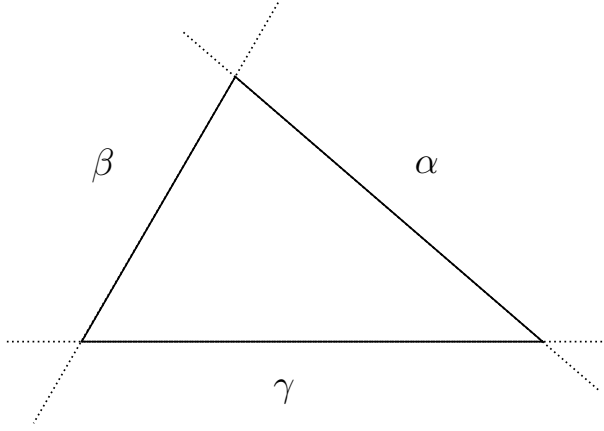

\centering
\begin{lpic}{triangle(8cm)}
\Large
\lbl[]{43,5;$\gamma$}
\lbl[]{65,40;$\alpha$}
\lbl[]{15,40;$\beta$}
\end{lpic}
\caption{A triangle in $X$ forms a complementary component of a Markov triple on $X$.}
\label{fig:triangle}
\end{figure}

\begin{proposition}
\label{prop:Dehn}
Suppose that $\alpha$, $\beta$, and $\gamma$ are a trio of simple closed geodesics on $X$ that pairwise intersect once, 
and suppose that a triangular complementary component sees $(\alpha,\beta,\gamma)$ in counterclockwise order, as in \autoref{fig:triangle}. Then the half-length function $\ell$ satisfies
\begin{align}
\label{eq:Dehn length}
\ell\left(\tau_\gamma^k(\beta)\right)
&= (k+1)
\ell(\gamma)+ \log\left( \frac32 m_{\alpha} - \frac{9m_\gamma m_\alpha-6m_\beta}{2\sqrt {9m_\gamma^2-4}}\right)+ O\left(e^{-2k\ell(\gamma)}\right)
\end{align}
as $k\to \infty$. (The implied constant in the last term depends on $\alpha$, $\beta$, and $\gamma$, but is independent of $k$.)
\end{proposition}

\begin{remark}
In terms of hyperbolic geometry, \autoref{prop:Dehn} is the computation of the asymptotic behavior of half-lengths under Dehn twisting. 
Unsurprisingly, the leading order term looks like $k$ times the half-length of $\gamma$. 
What is important for our purposes in this calculation is the constant term. 
\end{remark}

\begin{figure}[h]
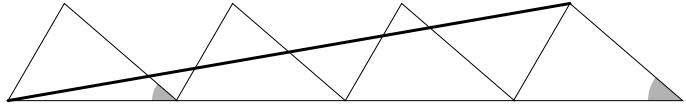

\centering
\begin{lpic}{lawCosines(9cm)}
\end{lpic}
\caption{The hyperbolic law of cosines can be applied to this triangle.}
\label{fig:cosines}
\end{figure}

\begin{proof}
We use the hyperbolic law of cosines on the triangle pictured in \autoref{fig:cosines}, with lengths $(k+1)\ell(\gamma)$, $\ell(\alpha)$, and $\ell(\tau_\gamma^k(\beta))$, the latter facing angle $\theta$:
\[
\cosh \ell(\tau_\gamma^k(\beta)) = \cosh \left((k+1)\ell(\gamma) \right)\cosh \ell(\alpha) - \sinh \left( (k+1)\ell(\gamma) \right) \sinh \ell(\alpha) \cos \theta~.
\]
Because $\cosh^{-1}\left( A_1 \cosh x+A_2 \sinh x\right) = x+ \log(A_1+A_2) + O\left( e^{-2x}\right)$ as $x\to\infty$,
\begin{align*}
\ell(\tau_\gamma^k(\beta))& =(k+1)\ell(\gamma) + \log(\cosh\ell(\alpha) - \sinh\ell(\alpha)\cos\theta)+O\left( e^{-2k\ell(\gamma)}\right)\\
& =(k+1)\ell(\gamma) + \log\left( \frac32 m_\alpha - \sinh \ell(\alpha) \cos \theta \right) + O\left( e^{-2k\ell(\gamma)}\right)~.
\end{align*}
Finally, applying the hyperbolic law of cosines to the triangle pictured in \autoref{fig:triangle}, we find that 
\[
\cosh \ell(\beta) = \cosh \ell(\alpha) \cosh \ell(\gamma) - \sinh\ell(\alpha) \sinh\ell(\gamma) \cos \theta~,
\]
whence $\sinh\ell(\alpha) \cos \theta = \frac{\cosh\ell(\alpha)\cosh \ell(\gamma) -\cosh \ell(\beta)}{\sinh\ell(\gamma)}$.
Substituting in and simplifying,
\begin{align*}
\ell(\tau_\gamma^k(\beta)) &= (k+1) \ell(\gamma) + \log\left( \frac32 m_\alpha - \frac{\cosh\ell(\alpha)\cosh \ell(\gamma) -\cosh \ell(\beta)}{\sinh\ell(\gamma)}\right) + O\left( e^{-2k\ell(\gamma)}\right) \\
&=(k+1) \ell(\gamma) + \log\left( \frac32 m_\alpha - \frac{\frac32m_\alpha \frac32 m_\gamma -\frac32 m_\beta}{\sqrt{ \frac94 m_\gamma^2 -1}}\right)+ O\left( e^{-2k\ell(\gamma)}\right) \\
&= (k+1) \ell(\gamma) + \log\left( \frac32 m_\alpha - \frac{9m_\alpha m_\gamma -6m_\beta}{2\sqrt{ 9 m_\gamma^2 -4}}\right)+ O\left( e^{-2k\ell(\gamma)}\right)~.\qedhere
\end{align*}
\end{proof}

\bigskip
\section{Derivative calculations}

Here we perform the computations for the proofs of \autoref{prop:psi derivatives} and \autoref{prop:corners}.
 
Let $\frac{r_1}{s_1}<\frac pq<\frac{r_2}{s_2}$ span a Farey triangle, with corresponding Markov triple $(n_1,n,n_2)$, and so that $n=m_{T(p/q)}$ and $n_i=m_{T(r_i/s_i)}$. (Recall that $T(p/q) = \frac p{p-q}$.)

The left-derivative of $\Psi(t)$ at $t=\frac pq$ may be computed as 
\[
\left.\frac{\upd }{\upd t}\right|_{t=\frac pq^-} \Psi(t) = \lim_{\frac uv \to \frac pq^-} \frac{\Psi\left( \frac pq\right) -\Psi\left( \frac uv \right)}{\frac pq - \frac uv} = \lim_{k\to \infty} \frac{ \Psi\left( \frac pq \right) - \Psi \left( \frac{u_k}{v_k} \right)}{\frac pq - \frac{u_k}{v_k}}~,
\]
where the last expression involves any sequence $\left(\frac {u_k}{v_k}\right)_{k=1}^\infty$ that approaches $\frac pq$ from the left. 
Since $\frac{r_1}{s_1}<\frac pq$, such a sequence can be obtained by Dehn twisting: let 
\[
\frac{u_k}{v_k} = \frac{r_1+kp}{s_1+kq}
\]
for each $k\in \bN$.
Now we have
\[
\frac pq - \frac{u_k}{v_k} = \frac pq - \frac{r_1+kp}{s_1+kq} = \frac 1{(s_1+kq)q} ~.
\]
In order to compute $\Psi\left( \frac{u_k}{v_k}\right)$, we need to understand the hyperbolic length $\|T(v_k,u_k)\|_s$.
Fortunately, we have in \eqref{eq:Dehn length} a calculation of the hyperbolic length of the curve represented by 
\[
T(v_k,u_k)=T((s_1,r_1)+k(q,p))=\tau^k_{T(q,p)}T(s_1,r_1)~.
\]
By \autoref{prop:Dehn} we find:
\begin{align*}
\frac{ \Psi\left( \frac pq \right) - \Psi \left( \frac{u_k}{v_k} \right)}{\frac pq - \frac{u_k}{v_k}} &=
(s_1+kq)q \left( \frac{ \|T(q,p)\|_s}{2q}- \frac{\left \|\tau^k_{T(q,p)}T(s_1,r_1)\right\|_s}{2(s_1+kq)}\right) \\
&=\frac12(s_1+kq)\|T(q,p)\|_s - \frac12q\left\|\tau^k_{T(q,p)}T(s_1,r_1)\right\|_s \\
&=\frac12 s_1 \|T(q,p)\|_s + \frac12kq\|T(q,p)\|_s 
\\ &\hspace{.5cm}- 
q \cdot
\left[ (k+1)\frac12\|T(q,p)\|_s + 
\log\left( \frac{3n_2}2 - 
\frac{9nn_2-6n_1}{2\sqrt{9n_2^2 -4}}\right)+ O\left(e^{-k\|(q,p)\|_s}\right) \right]
\\
&=(s_1-q)\cosh^{-1}\left(\frac{3n}2 \right)
-q\log\left( \frac{3n_2}2 - 
\frac{9nn_2-6n_1}{2\sqrt{9n_2^2 -4}}\right)+ O\left(e^{-k\|(q,p)\|_s}\right)\\
&= -s_2\cosh^{-1}\left(\frac{3n}2 \right)
-q\log\left( \frac{3n_2}2 - 
\frac{9nn_2-6n_1}{2\sqrt{9n_2^2 -4}}\right)+ O\left(e^{-k\|(q,p)\|_s}\right) 
\end{align*}
The limit is evident as $k\to \infty$, and the computation is complete.

Of course, the exact same computation can be carried out for the right-derivative as well. We leave the details to the reader.

The transformation $f:(x,y)\mapsto\left(\frac{1-x}{2y},\frac x{2y}\right)$ maps the graph of $\Psi$ to the portion of the unit sphere, with respect to the stable norm, that lies inside $\{(q,p):q\ge p\ge 0\}$. Indeed,
\begin{align*}
\| f(t,\Psi(t)) \|_s &=\left\|\left( \frac{1-t}{2\Psi(t)} , \frac t{2\Psi(t)} \right) \right\|_s = \frac1{2\Psi(t)}\|(1-t,t)\|_s~.
\end{align*}
When $t=p/q\in\bQ$, by \eqref{eq:psi-snorm} we find that
\[
\frac1{2\Psi(t)} \|(1-t,t)\|_s = \frac1{2\Psi(p/q)} \left\|\left(1-\frac pq, \frac pq \right)\right\|_s = \frac1{\Psi(p/q)} \cdot \frac{\|T(q,p)\|_s}{2q} = 1~.
\]
By continuity, $\|f(t,\Psi(t))\|_s=1$ for all $t\in[0,1/2]$.

It is routine to compute the derivative $df$ at $(t,\Psi(t))$, and \autoref{prop:corners} follows easily.

\bigskip

\section{Proof of \autoref{thm:main thm}}

We show how \autoref{thm:main thm} can be deduced from \autoref{thm:fixed point}.

\begin{proof}[Proof of \autoref{thm:main thm}]
Let $(q,p),(q',p')\in\cQ$ and $r=\frac{q'-q}{p'-p}$, and suppose without loss of generality that $q'< q$.
Because $\sigma_-=\mu^+_{0/1}$ and $\sigma_+=\mu^-_{1/1}$ (see \autoref{ex:minus} and \autoref{ex:plus}), convexity of the stable norm sphere implies that $\mu^-_{p/q},\mu^+_{p/q},\mu^-_{p'/q'},\mu^+_{p'/q'}\in[\sigma_-,\sigma_+]$. 
If $r<\sigma_-\le \mu^+_{p'/q'}$ then \autoref{thm:fixed point} implies that $(q,p)\prec_\cM (q',p')$, so $\prec_\cM$ is decreasing with $q$, demonstrating \autoref{it:sigma-}.
If $r>\sigma_+\ge\mu^-_{p'/q'}$ then \autoref{thm:fixed point} implies that $(q',p')\prec_\cM (q,p)$, so $\prec_\cM$ is increasing with $q$, demonstrating \autoref{it:sigma+}.

Now suppose that $r\in(\sigma_-,\sigma_+)$, and let $\cL$ be a line of slope $r$. Let $P$ and $Q$, respectively, be the intersections of $\cL$ with the lines through the origin of slopes $0$ and $1$. Because $r>\sigma_- = \mu^+_{0/1}$, near $P$ the line $\cL$ is inside the ball of radius $ \|P\|_s$ about the origin; because $r<\sigma_+=\mu^-_{1/1}$, near $Q$ it is inside the ball of radius $ \|Q\|_s$ about the origin. It follows that, restricted to the line $\cL$, as a function of 
$q$ the norm $\|{\cdot}\|_s$ is decreasing near $Q$ and increasing near $P$. 
By convexity, there is a point $O=(q_o,p_o)$ on the line segment $QP$ so that, as a function of $q$, $\|{\cdot}\|_s$ is decreasing on $QO$ and increasing on $OP$.

Recall that $\cL_t$ indicates $t\cdot \cL$, the scaling of $\cL$ by $t$.
It remains to explore $\prec_\cM$ on $\cL_t$, or, equivalently, $\|{\cdot}\|_s$ on $\cL_t \cap \cQ$.
Let $\cQ'=\{(q,p)\in\bZ^2:q\ge p\ge 0\}$.
Suppose that $\cL_t\cap \cQ' = \{(q_1,p_1),\ldots,(q_n,p_n)\}$ with $q_1<\ldots<q_n$.
Strict convexity of $\|{\cdot}\|_s$ implies that there is some $j\in\{1,\ldots,n\}$ so that 
\begin{equation}
\label{eq:antimodal}
\|(q_1,p_1)\|_s > \ldots > \|(q_j,p_j)\|_s < \ldots < \|(q_n,p_n)\|_s~.
\end{equation}
Moreover, for $t$ large enough \eqref{eq:antimodal} holds with $2\le j\le n-1$:
The points $\left\{\frac 1t (q_i,p_i) \right\} $ are dense in $PQ$, so by continuity of $\|{\cdot}\|_s$ there is a $t$ large enough so that there are points $\frac1t(q_i,p_i)$ in both $QO$ and $OP$, i.e.~$\|\frac1t(q_1,p_1)\|_s> \|\frac1t(q_2,p_2)\|_s$ and $\|\frac1t(q_{n-1},p_{n-1})\|_s<\|\frac1t(q_n,p_n)\|_s$.
This demonstrates non-monotonicity with respect to $q$ of $\|{\cdot}\|_s$ on $\cL_t\cap \cQ'$ (and, in fact, \emph{strict} antimodality of $\|{\cdot}\|_s$ on $\cL_t\cap\cQ'$).

In order to deduce strict antimodality (and non-monotonicity) for $\prec_\cM$, it remains to go from $\cL_t\cap \cQ'$ to $\cL_t\cap\cQ$, i.e.~one needs choices for $t$ so that there are many relatively prime points $(q_i,p_i)\in \cL_t\cap\cQ$. 

The following demonstrates that such choices exist for arbitrarily large $t$:
Choose $t$ and $k$ so that $(k,k-1)\in \cL_t\cap\cQ$, for a large $k\in\bN$ relatively prime to $v$, and let 
\begin{equation}
\label{eq:pts}
(q_j,p_j)=(k,k-1)+(j-1)(v,-u) \text{\ \ for \ }j=1,\ldots, K ~,
\end{equation} 
where $K=\left\lfloor  \frac{k-1}u \right\rfloor$.
The Siegel-Walfisz Theorem provides an asymptotic count of $\pi(x;v,k)$, the number of primes $\le x$ which are $\equiv k(\mathrm{mod}\, v)$. 
Provided $x$ is large relative to $v$ (indeed, for us $v$ is fixed while $x\to \infty$) and $h$ is large relative to $x$ (for us, $h\gtrsim x$), one has
\[
\pi(x;v,k) - \pi(x-h;v,k) \sim \frac h{\phi(v) \log x}~.
\]  
(See \cite[Cor.~11.19]{MV-book} for a precise statement of the Siegel-Walfisz Theorem.) 

We find that for any $\delta>0$  and $h=\delta x$ there exists $N$ so that
\[
\pi((1+\delta)x;v,k) - \pi(x;v,k) \ge 2
\]
for all $x\ge N$. 
As $t$ and $k$ satisfy a linear relationship (that is $C_1 k \le t\le C_2 k$ for some constants $C_1$ and $C_2$),
we may apply this estimate to both the intervals $[k,tq_o]$ and $[tq_o, q_K]$, and one obtains several points with prime $q$-coordinate inside $t\cdot QO$ and $t\cdot OP$. 
This demonstrates that the points $\{(q_i,p_i)\}$ defined in \eqref{eq:pts} contain points from $\cQ$ for which $\prec_\cM$ both decreases and increases with respect to $q$. 
Hence we have strict antimodality with respect to $q$ for $\cL_t\cap \cQ$. \qedhere

\end{proof}

\begin{remark}
\label{rem:strict antimodality}
The reader may observe that \eqref{eq:antimodal} implies that there exists $T>0$ so that, for all $t\ge T$, the stable norm $\|{\cdot}\|_s$ is strictly antimodal with respect to $q$ on the intersection $\cL_t\cap \cQ'$.
\end{remark}

\begin{remark}
\label{rem:Conj612}
Lee-Li-Rabideau-Schiffler claim 
more than the statement of \autoref{thm:main thm}\ref{it:nope}.
Namely, they extend $\prec_\cM$ from $\cQ$ to $\cQ'$ (see also \cite{LPTV}). The latter is induced by something they call the `Markov distance' -- though one should note that the Markov distance is not a metric, as it does not obey the triangle inequality \cite[Rem.\,3.7]{LLRS}. 
With the notation that $m_{q,p}$ stands for the Markov distance from $(0,0)$ to $(q,p)$, one has $m_{q,p}=m_{p/q}$ when $\gcd(q,p)=1$. 
%
A formula computing $m_{q,p}$ in terms of $m_{q',p'}$, where $(q',p')=\frac1{\gcd(q,p)}(q,p)$, can be found in \cite[Lem.\,6.2]{LLRS}. 
One finds that
\begin{equation}
\label{eq:Markov distance}
m_{q,p}=\frac1{L_{q',p'}} \cdot 2\sinh\left(\frac12\|(q,p)\|_s\right)~,
\footnote{This formula corrects the false claim that $m_{q,p}=\frac23\cosh\frac12\|(q,p)\|_s$ found in a previous version of this manuscript. The previous `Corollary 4.2', linking Markov distances to traces of non-primitive curves, is likewise withdrawn.}
\end{equation}
where $L_{s,r}=\sqrt{ 9-\frac4{m^2}}$ for $m=m_{r/s}$ (that is, $L_{s,r}$ is the \emph{Lagrange number} associated to Markov number $m_{r/s}$ -- the Lagrange numbers are fundamental invariants in Diophantine approximation \cite{Aigner, CF}).

In \cite[Conj.\,6.12]{LLRS}, it is conjectured that there is a $T$ so that, for all $t\ge T$, the Markov distance is strictly antimodal with $q$ on the intersection $\cL_t\cap \cQ'$. 
As $\sinh$ is an increasing function, it is tempting to prove this claim by combining \autoref{rem:strict antimodality} and \eqref{eq:Markov distance}. 
However, the Lagrange number $L_{q',p'}$ provides a significant obstacle. 
Indeed, we point out that the orders on $\cQ'$ induced by the Markov distance (in which $(q,p)<(r,s)$ whenever $m_{q,p}<m_{r,s}$) and by the stable norm $\|{\cdot}\|_s$ (in which $(q,p)<(r,s)$ whenever $\|(q,p)\|_s<\|(r,s)\|_s$) do \textbf{not} coincide in general: one has
\[
m_{9,0}=2584 > 2378 = m_{5,5} \ \ \text{ while } \ \ \|(9,0)\|_s=8.66\ldots < 8.81\ldots = \|(5,5)\|_s~.
\]
This is in stark contrast to the situation on $\cQ$, in which $m_{q,p}$ and $\|(q,p)\|_s$ are related by an increasing function \eqref{eq:length-Markov}.
Because the orders induced by the Markov distance and the stable norm do not coincide, the proof above of strict antimodality of the stable norm on $\cL_t\cap \cQ'$ for $t$ sufficiently large does not translate to the same for the Markov distance.
Thus our proof stops short of demonstrating the full statement of \cite[Conj.\,6.12]{LLRS}.
\end{remark}

\bigskip

\bibliographystyle{alpha}

\begin{thebibliography}{LPTV21}

\bibitem[Aig15]{Aigner}
Martin Aigner.
\newblock {\em Markov's theorem and 100 years of the uniqueness conjecture}.
\newblock Springer, 2015.

\bibitem[BBI01]{BBI}
Dmitri Burago, Yuri Burago, and Sergei Ivanov.
\newblock {\em A course in metric geometry}, volume~33.
\newblock American Mathematical Soc., 2001.

\bibitem[CF89]{CF}
Thomas~W Cusick and Mary~E Flahive.
\newblock {\em The {M}arkoff and {L}agrange spectra}.
\newblock Number~30. American Mathematical Soc., 1989.

\bibitem[Coh55]{Cohn}
Harvey Cohn.
\newblock Approach to {M}arkoff's minimal forms through modular functions.
\newblock {\em Annals of Mathematics}, pages 1--12, 1955.

\bibitem[Cul86]{Culler}
Marc Culler.
\newblock Lifting representations to covering groups.
\newblock {\em Advances in Mathematics}, 59(1):64--70, 1986.

\bibitem[EMM19]{EMM}
Alex Eskin, Maryam Mirzakhani, and Amir Mohammadi.
\newblock Effective counting of simple closed geodesics on hyperbolic surfaces.
\newblock arXiv:1905.04435, 2019.

\bibitem[Foc97]{Fock}
Vladimir Fock.
\newblock Dual {T}eichm\" uller spaces.
\newblock arXiv:dg-ga/9702018, 1997.

\bibitem[Fro13]{Frobenius}
Georg~Ferdinand Frobenius.
\newblock {\em \"{U}ber die {M}arkoffschen zahlen}.
\newblock K\"onigliche Akademie der Wissenschaften, 1913.

\bibitem[Gol03]{Goldman}
William~M Goldman.
\newblock The modular group action on real {SL}(2)--characters of a one-holed
  torus.
\newblock {\em Geometry \& Topology}, 7(1):443--486, 2003.

\bibitem[Gor81]{Gorshkov}
DS~Gorshkov.
\newblock Geometry of {L}obachevskii in connection with certain questions of
  arithmetic.
\newblock {\em Journal of Soviet Mathematics}, 16(1):788--820, 1981.


\bibitem[Hin20]{Hines}
Robert Hines.
\newblock An infinite product on the {T}eichm\"uller space of the once-punctured torus.
\newblock arXiv:2001.05557, 2020.


\bibitem[LLRS20]{LLRS}
Kyungyong Lee, Li~Li, Michelle Rabideau, and Ralf Schiffler.
\newblock On the ordering of the {M}arkov numbers.
\newblock arXiv:2010.13010, 2020.

\bibitem[LPTV21]{LPTV}
Cl{\'e}ment Lagisquet, Edita Pelantov{\'a}, S{\'e}bastien Tavenas, and Laurent
  Vuillon.
\newblock On the {M}arkov numbers: fixed numerator, denominator, and sum
  conjectures.
\newblock {\em Advances in Applied Mathematics}, 130:102227, 2021.

\bibitem[Mar79]{Markoff}
Andrey Markoff.
\newblock Sur les formes quadratiques binaires ind{\'e}finies.
\newblock {\em Mathematische Annalen}, 15(3):381--406, 1879.

\bibitem[McS04]{McShane-Weierstrass}
Greg McShane.
\newblock Weierstrass points and simple geodesics.
\newblock {\em Bulletin of the London Mathematical Society}, 36(2):181--187,
  2004.

\bibitem[McS21]{McShane}
Greg McShane.
\newblock Convexity and {A}igner's conjectures.
\newblock arXiv:2101.03316, 2021.

\bibitem[MR95a]{MR1}
Greg McShane and Igor Rivin.
\newblock A norm on homology of surfaces and counting simple geodesics.
\newblock {\em International Mathematics Research Notices}, 1995(2):61--69,
  1995.

\bibitem[MR95b]{MR2}
Greg McShane and Igor Rivin.
\newblock Simple curves on hyperbolic tori.
\newblock {\em Comptes rendus de l'Acad{\'e}mie des sciences. S{\'e}rie 1,
  Math{\'e}matique}, 320(12):1523--1528, 1995.

\bibitem[MV07]{MV-book}
Hugh~L Montgomery and Robert~C Vaughan.
\newblock {\em Multiplicative number theory I: Classical theory}.
\newblock Number~97. Cambridge university press, 2007.

\bibitem[Pet06]{Petrov}
Fedor~Vladimirovich Petrov.
\newblock On the number of rational points on a strictly convex curve.
\newblock {\em Functional Analysis and Its Applications}, 40(1):24--33, 2006.

\bibitem[Pla99]{Plagne}
Alain Plagne.
\newblock A uniform version of {J}arn{\i}k's theorem.
\newblock {\em Acta Arith}, 87:255--267, 1999.

\bibitem[Rab18]{Rabideau}
Michelle Rabideau.
\newblock {\em Continued Fractions in Cluster Algebras, Lattice Paths and
  {M}arkov Numbers}.
\newblock PhD thesis, U.~of Conn., 2018.

\bibitem[RS20]{RS}
Michelle Rabideau and Ralf Schiffler.
\newblock Continued fractions and orderings on the {M}arkov numbers.
\newblock {\em Advances in Mathematics}, 370:107231, 2020.

\bibitem[SD74]{Swinnerton-Dyer}
HPF Swinnerton-Dyer.
\newblock The number of lattice points on a convex curve.
\newblock {\em Journal of Number Theory}, 6(2):128--135, 1974.

\bibitem[Ser85]{Series}
Caroline Series.
\newblock The geometry of {M}arkoff numbers.
\newblock {\em The {M}athematical {I}ntelligencer}, 7(3):20--29, 1985.

\bibitem[SV17]{SpV}
Kathryn Spalding and AP~Veselov.
\newblock Lyapunov spectrum of {M}arkov and {E}uclid trees.
\newblock {\em Nonlinearity}, 30(12):4428, 2017.

\bibitem[SV19]{SoV}
Alfonso Sorrentino and Alexander~P Veselov.
\newblock Markov numbers, {M}ather's $\beta$ function and stable norm.
\newblock {\em Nonlinearity}, 32(6):2147, 2019.

\bibitem[Zag82]{Zagier}
Don Zagier.
\newblock On the number of {M}arkoff numbers below a given bound.
\newblock {\em Mathematics of Computation}, 39(160):709--723, 1982.

\end{thebibliography}

\bigskip

\end{document}